\documentclass[smallextended,referee,envcountsect,]{svjour3}
\smartqed
\usepackage{graphicx}

\usepackage{amsmath}
\usepackage{amssymb}
\usepackage{amsfonts}
\usepackage{mathtools}
\usepackage{bm}
\usepackage{adjustbox}

\usepackage{empheq}
\usepackage{booktabs}
\usepackage{tabularx}

\usepackage{algorithm}
\usepackage{algpseudocode}

\algrenewcommand\algorithmicrequire{\textbf{Input:}}

\usepackage{graphicx}

\usepackage{hyperref}

\DeclareMathOperator{\Exp}{Exp}
\DeclareMathOperator{\Log}{Log}

\spnewtheorem{assumption}{Assumption}[section]
    {\bfseries}{\itshape}

\smartqed
\smartqed

\journalname{JOTA}

\begin{document}

\title{Region-Preserving Accelerated Gradient Methods for Geodesically Convex Optimization}

\author{Derun Zhou}

\institute{Derun Zhou \at
National Institute of Informatics (NII)\\
2-1-2 Hitotsubashi, Chiyoda-ku,
Tokyo 101-8430, Japan\\
The Graduate University for Advanced Studies
(SOKENDAI)\\
Hayama, Kanagawa 240-0193, Japan\\
\email{zhouderun@nii.ac.jp}
}

\date{Received: date / Accepted: date}

\maketitle

\begin{abstract}
Accelerated Riemannian optimization methods often require
a priori assumptions that their iterates remain in a region
where geometric bounds hold. We develop a parameterized
inertial dynamical system and an explicit accelerated
gradient method for smooth geodesically convex optimization
on Riemannian manifolds with locally bounded sectional
curvature. Our setting allows both positive and negative
curvature without requiring completeness of the ambient
manifold. Starting from a curvature-corrected continuous
model, we construct the discrete updates by eliminating
indefinite cross terms in the Lyapunov energy difference
and derive verifiable conditions for a fixed step size.
The resulting energy analysis establishes preservation
of the initial geodesic ball without assuming iterate
containment in advance. We further establish accelerated
convergence rates for both geodesically convex and
strongly geodesically convex objectives, without requiring
restarts. 
\end{abstract}




\section{Introduction}
\label{sec:introduction}

Riemannian optimization concerns the minimization of
functions defined on smooth manifolds, where the
underlying geometric structure can be exploited to
develop efficient optimization algorithms. Of
particular interest are geodesically convex problems,
which arise in covariance estimation, matrix
optimization, and the computation of Riemannian
means. We consider
\begin{equation}
\label{eq:intro-problem}
\min_{x\in\mathcal N\subset\mathcal M} f(x),
\end{equation}
where $\mathcal M$ is a finite-dimensional Riemannian
manifold, $\mathcal N$ is a bounded, uniquely
geodesically convex domain, and $f$ is geodesically
convex or strongly geodesically convex.

In Euclidean optimization, Nesterov's accelerated
gradient methods improve the convergence rate of
gradient descent from $O(k^{-1})$ to $O(k^{-2})$
for smooth convex objectives, and achieve accelerated
linear convergence for strongly convex objectives.
Extending these guarantees to Riemannian manifolds
is challenging because the nonlinear geometry
introduces distortion between tangent spaces.
Consequently, the Euclidean estimate-sequence and
Lyapunov arguments cannot be transferred directly
to the Riemannian setting.

Several accelerated Riemannian methods have been
developed to address this difficulty.
\cite{ZhangSra2018} established accelerated
convergence for smooth, strongly geodesically
convex objectives when the initial point lies
in a sufficiently small neighborhood of the
minimizer. \cite{AhnSra2020} subsequently
developed a Riemannian acceleration framework
based on estimate sequences and metric distortion.
For geodesically convex objectives,
\cite{alimisis2021momentum} proposed a
momentum-based method with curvature-dependent
convergence guarantees. However, its general
bound contains an additional $O(k^{-1})$ term.

A major development was made by
\cite{KimYang2022}, who proposed the Riemannian
Nesterov accelerated gradient methods RNAG-C and
RNAG-SC. Under bounded sectional curvature and
a bounded geodesically convex domain, these
methods attain the Euclidean accelerated
convergence orders for both geodesically convex
and strongly convex objectives. Nevertheless,
their convergence analysis assumes a priori
that all iterates $x_k$ and intermediate points
$y_k$ remain within the prescribed domain.
Without an independent mechanism ensuring this
property, the local curvature and smoothness
bounds used in the analysis may cease to apply.

The iterate-containment issue was first addressed
under constant sectional curvature by
\cite{MartinezRubio2022}, who established full
acceleration in hyperbolic and spherical spaces.
Their method transforms the problem into
constrained Euclidean optimization via
a geodesic map, restricting the iterates to
$\mathcal B_R=\Exp_{x_0}(\overline B(0,R))$,
where $R\ge d(x_0,x^\star)$.
For positive curvature $K>0$, the radius must
additionally satisfy $R<\pi/(2\sqrt K)$.
However, this approach relies on constant
sectional curvature and constrained optimization.

Subsequently, \cite{MartinezRubioPokutta2023}
developed an accelerated Riemannian inexact
proximal-point framework using constrained
proximal subproblems and tangent-space projections
to control the iterates. For smooth objectives,
the proximal subproblems are approximated by
first-order methods on bounded geodesic balls.
Their first-order acceleration theorem applies
to Hadamard manifolds satisfying
$\nabla\mathrm{Rm}=0$ and requires smoothness
and geodesic convexity on
$\overline B_{\mathcal M}(x^\star,3R)$,
where $R\ge d(x_0,x^\star)$.
The strongly convex complexity analysis also
employs a restart reduction.
While both approaches guarantee iterate
containment, they rely on constrained
computations and additional geometric assumptions.

These developments motivate the question of
whether an explicit Riemannian accelerated
gradient method can preserve a prescribed
region without assuming iterate containment
a priori or introducing constrained proximal
subproblems. In this paper, we address this
question through a curvature-corrected
continuous-time model and an associated
Lyapunov-based discretization.

Our analysis adopts the bounded-region
geometric setting of \cite{KimYang2022},
but replaces the a priori iterate-containment
assumption with
\begin{equation}
\label{eq:intro-initial-ball}
\overline B_{\mathcal M}(x^\star,r_0)
\Subset\mathcal N,
\qquad r_0=d(x_0,x^\star).
\end{equation}
Here $\Subset$ denotes compact containment.
Unlike the Hadamard setting, our assumptions
allow sectional curvature of either sign and
do not require the ambient manifold to be
complete. Under these conditions, we establish
that the continuous trajectory, the discrete
iterates, and the entire geodesic update
segments remain within the initial ball.

Our approach is motivated by continuous-time
models of accelerated optimization
\cite{Alimisis2020,KimYang2023Kernels}
and by the Lyapunov-based discretization
of ITEM dynamics \cite{Ushiyama2025}.
We introduce a parameterized Riemannian ITEM
ODE whose damping incorporates curvature
information and whose Lyapunov function
controls both the objective residual and
the distance to the minimizer.
Rather than directly discretizing the ODE,
we determine the discrete update coefficients
by canceling the mixed terms in the difference
of the corresponding Lyapunov function.
This leads to an explicit accelerated gradient
method with one new gradient evaluation per
iteration, without projections, proximal
subproblems, or restarts.

\noindent\textbf{Contributions.}
The main contributions of this paper are summarized
as follows.

\begin{enumerate}
\item We introduce a curvature-corrected,
$\sigma$-parameterized Riemannian ITEM ODE and
construct a Lyapunov function that establishes
energy dissipation, initial-ball invariance,
and accelerated continuous-time convergence
for geodesically convex and strongly convex
objectives.

\item We derive an explicit Riemannian
accelerated gradient method from the discrete
Lyapunov energy. The update coefficients are
determined by canceling the mixed terms in
the energy difference, yielding a unified
scheme for the convex and strongly convex
settings without restarts or proximal
subproblems.

\item We develop a uniform matrix-based
step-size condition that guarantees the
decrease of the discrete Lyapunov function.
We establish the existence of a positive
certified fixed step size under admissible
curvature-dependent parameter choices,
including the Euclidean choice $s=1/L$.
The certification is performed before
the iterations and does not require
step-by-step adjustment.

\item We establish initial-ball invariance
for the entire discrete update geodesics,
without assuming a priori that the iterates
remain in the working domain. Consequently,
the method attains an $O(k^{-2})$ function-value
rate for geodesically convex objectives and
an accelerated linear rate for strongly
geodesically convex objectives, under bounded
sectional curvature of either sign.
\end{enumerate}

The remainder of the paper is organized as
follows. Section~\ref{sec:preliminaries}
reviews the necessary Riemannian preliminaries.
Section~\ref{sec:assumptions} introduces
the geometric, objective, and initial-region
assumptions.
Section~\ref{sec:continuous-model} develops
the continuous-time dynamics and Lyapunov analysis.
Section~\ref{sec:sigma-item-discretization}
presents the discretization, fixed-step
certification, and convergence results.
The final sections report numerical experiments
and concluding remarks.

\section{Preliminaries on Riemannian Geometry}
\label{sec:preliminaries}

We briefly recall the basic notions of Riemannian geometry and
geodesic convexity used throughout the paper.

Let $(\mathcal M,g)$ be a finite-dimensional Riemannian manifold.
For each $x\in\mathcal M$, let $T_x\mathcal M$ denote the tangent space
at $x$. The Riemannian metric induces the inner product and norm
\[
\langle u,v\rangle_x:=g_x(u,v),
\qquad
\|v\|_x:=\sqrt{\langle v,v\rangle_x},
\qquad
u,v\in T_x\mathcal M,
\]
and the tangent bundle is
$T\mathcal M:=\bigsqcup_{x\in\mathcal M}T_x\mathcal M$.

For a piecewise smooth curve $\gamma:[a,b]\to\mathcal M$, its length,
the Riemannian distance, and the diameter of a nonempty set
$\mathcal N\subseteq\mathcal M$ are defined, respectively, by
\[
\ell(\gamma):=\int_a^b\|\dot\gamma(t)\|_{\gamma(t)}\,dt,
\qquad
d(x,y):=\inf_{\gamma:x\to y}\ell(\gamma),
\qquad
\operatorname{diam}(\mathcal N):=
\sup_{x,y\in\mathcal N}d(x,y).
\]

Let $\nabla$ denote the Levi--Civita connection. A smooth curve
$\gamma$ is a geodesic if $\nabla_{\dot\gamma}\dot\gamma=0$.
For $x\in\mathcal M$ and $v\in T_x\mathcal M$, the exponential map is
defined by
\[
\Exp_x(v):=\gamma_v(1),
\]
where $\gamma_v$ is the geodesic satisfying
$\gamma_v(0)=x$ and $\dot\gamma_v(0)=v$. Whenever $\Exp_x$ is
invertible on the relevant domain, its inverse is denoted by $\Log_x$.
If $x$ and $y$ are joined by a unique minimizing geodesic, we write
\[
\gamma_{x,y}(t)
:=
\Exp_x\bigl(t\Log_x y\bigr),
\qquad
t\in[0,1],
\qquad
d(x,y)=\|\Log_x y\|_x.
\]

For such $x$ and $y$, let
$P_{x\to y}:T_x\mathcal M\to T_y\mathcal M$ denote parallel transport
along the unique minimizing geodesic. It is a linear isometry satisfying
\[
\langle P_{x\to y}u,P_{x\to y}v\rangle_y
=
\langle u,v\rangle_x,
\qquad
\|P_{x\to y}v\|_y=\|v\|_x,
\qquad
P_{y\to x}=P_{x\to y}^{-1}.
\]

For a differentiable function $f:\mathcal M\to\mathbb R$, its
Riemannian gradient $\operatorname{grad}f(x)\in T_x\mathcal M$ is
defined by
\[
Df(x)[v]
=
\langle\operatorname{grad}f(x),v\rangle_x,
\qquad
v\in T_x\mathcal M.
\]
If $f$ is twice differentiable, its Riemannian Hessian is
\[
\operatorname{Hess}f(x)[v]
:=
\nabla_v\operatorname{grad}f,
\]
and, for any smooth curve $\gamma$,
\[
\frac{d}{dt}f(\gamma(t))
=
\left\langle
\operatorname{grad}f(\gamma(t)),\dot\gamma(t)
\right\rangle_{\gamma(t)}.
\]
For linearly independent $u,v\in T_x\mathcal M$,
$\operatorname{Sec}_x(u,v)$ denotes the sectional curvature of the
two-dimensional subspace spanned by $u$ and $v$; when no ambiguity
arises, we simply write $\operatorname{Sec}$.

A set $\mathcal N\subseteq\mathcal M$ is called
\emph{geodesically convex} if every $x,y\in\mathcal N$ can be joined
by a minimizing geodesic contained in $\mathcal N$, and
\emph{geodesically uniquely convex} if this minimizing geodesic is
unique.

A function $f:\mathcal N\to\mathbb R$ is called
\emph{geodesically convex} if
\[
f\bigl(\gamma_{x,y}(t)\bigr)
\leq
(1-t)f(x)+tf(y),
\qquad
x,y\in\mathcal N,\quad t\in[0,1].
\]
For differentiable $f$, this is equivalently characterized by
\begin{equation}
\label{eq:geodesic-convexity}
f(y)
\geq
f(x)
+
\left\langle
\operatorname{grad}f(x),\Log_x y
\right\rangle_x,
\qquad
x,y\in\mathcal N.
\end{equation}

More generally, $f$ is
\emph{$\mu$-strongly geodesically convex}, with $\mu\geq0$, if
\begin{equation}
\label{eq:strong-geodesic-convexity}
f(y)
\geq
f(x)
+
\left\langle
\operatorname{grad}f(x),\Log_x y
\right\rangle_x
+
\frac{\mu}{2}d^2(x,y),
\qquad
x,y\in\mathcal N.
\end{equation}
The case $\mu=0$ reduces to geodesic convexity.

A differentiable function $f:\mathcal N\to\mathbb R$ is called
\emph{$L$-smooth} if its Riemannian gradient is $L$-Lipschitz with
respect to parallel transport:
\begin{equation}
\label{eq:riemannian-L-smoothness}
\left\|
\mathcal P_{x\to y}\operatorname{grad}f(x)
-
\operatorname{grad}f(y)
\right\|_y
\leq
L\,d(x,y),
\qquad
x,y\in\mathcal N.
\end{equation}
This condition implies the Riemannian descent lemma
\begin{equation}
\label{eq:riemannian-descent-lemma}
f(y)
\leq
f(x)
+
\left\langle
\operatorname{grad}f(x),\Log_x y
\right\rangle_x
+
\frac{L}{2}d^2(x,y),
\qquad
x,y\in\mathcal N.
\end{equation}

\section{Assumptions}
\label{sec:assumptions}
We work on a geodesically convex normal domain and impose
curvature bounds only on the initial geodesic ball, rather than
on the entire ambient manifold or a bounded-diameter working domain.
Unlike analyses that assume the iterates remain in the working
domain a priori \cite{alimisis2021momentum,KimYang2022}, we
establish this property from the initial-region condition below.

\begin{assumption}
\label{ass:geometry}
Let $(\mathcal M,g)$ be a finite-dimensional Riemannian manifold,
not necessarily complete. The set $\mathcal N\subset\mathcal M$
is open and uniquely geodesically convex: any two points of
$\mathcal N$ are joined by a unique minimizing geodesic contained
in $\mathcal N$. 
\end{assumption}

\begin{remark}
Under Assumption~\ref{ass:geometry}, for every
$x\in\mathcal N$, the restriction of the exponential
map
\[
\Exp_x:\Log_x(\mathcal N)\to\mathcal N
\]
is a diffeomorphism
\cite{KimYang2022,alimisis2021momentum}.
Global invertibility of $\Exp_x$ on $T_x\mathcal M$
is not required.
\end{remark}

\begin{assumption}
\label{ass:objective}
Let $f\in C^2(\mathcal N)$ be geodesically $L$-smooth,
with $L>0$, and $\mu$-strongly geodesically convex for
$0\le\mu\le L$ in the sense of
\eqref{eq:strong-geodesic-convexity}, where $\mu=0$
corresponds to geodesic convexity. Assume that $f$ admits
a minimizer $x^\star\in\mathcal N$, and write
$f^\star=f(x^\star)$.
\end{assumption}

\begin{assumption}
\label{ass:initial-region}
Let $x_0\in\mathcal N$ and set
$r_0:=d_{\mathcal M}(x_0,x^\star)$. Assume that
\begin{equation}
\label{eq:initial-region}
C_0:=\overline B_{\mathcal M}(x^\star,r_0)
\Subset\mathcal N,
\end{equation}
where $\overline B_{\mathcal M}(x^\star,r_0)
=\{x\in\mathcal M:d_{\mathcal M}(x,x^\star)\le r_0\}$
is the closed geodesic ball and $\Subset$ denotes compact
containment. 

We assume that the sectional curvature on $C_0$
is bounded by finite constants $K_{\min}\le K_{\max}$.
Specifically, for every $x\in C_0$ and every
two-dimensional subspace
$\Pi=\operatorname{span}\{u,v\}\subset T_x\mathcal M$,
\begin{equation}
\label{eq:sectional-curvature-bound}
K_{\min}
\le \operatorname{Sec}_x(\Pi)
:=\frac{\langle\mathrm{Rm}_x(u,v)v,u\rangle_x}
{\|u\|_x^2\|v\|_x^2-\langle u,v\rangle_x^2}
\le K_{\max},
\end{equation}
where $u,v$ are linearly independent and
$\mathrm{Rm}$ denotes the Riemannian curvature tensor.
If $K_{\max}>0$, we additionally assume
$r_0<\pi/\sqrt{K_{\max}}$.

The associated radial comparison constants are
\begin{equation}
\label{eq:curvature-comparison-constants}
\begin{aligned}
\delta_{r_0}&:=
\begin{cases}
1,&K_{\max}\le0,\\[1mm]
\sqrt{K_{\max}}r_0\cot(\sqrt{K_{\max}}r_0),
&K_{\max}>0,
\end{cases}\\[2mm]
\zeta_{r_0}&:=
\begin{cases}
\sqrt{-K_{\min}}r_0\coth(\sqrt{-K_{\min}}r_0),
&K_{\min}<0,\\[1mm]
1,&K_{\min}\ge0.
\end{cases}
\end{aligned}
\end{equation}
For $r_0=0$ both comparison constants are defined as $1$.
In all cases $\delta_{r_0}\le1\le\zeta_{r_0}$.
\end{assumption}

\begin{remark}
\label{rem:initial-region}
The curvature bounds in
Assumption~\ref{ass:initial-region} are imposed only
on $C_0$, not on the whole of $\mathcal N$ or $\mathcal M$.
Our analysis does not require
$\operatorname{diam}(\mathcal N)<\infty$ or an a priori
assumption that $x_k,y_k\in\mathcal N$. Instead, we prove
that the continuous trajectory and both complete update
geodesics remain in $C_0$. The compactness in
\eqref{eq:initial-region} is relevant when
$\mathcal M$ is not complete.
\end{remark}

\section{A parameterized Riemannian ITEM ODE}
\label{sec:continuous-model}

In this section, we introduce a parameterized Riemannian extension
of the ITEM ODE considered in \cite{KimYang2023Kernels} and study its
Lyapunov dissipation, region invariance, and convergence rates.
Its energy retains the negative squared-distance term of
\cite{Ushiyama2025}. We first introduce the 
Hessian comparison estimate, using only the curvature bounds on $C_0$.

\begin{lemma}[\cite{Alimisis2020}, Lemma~2]
\label{lem:hessian-comparison-squared-distance}
Under Assumptions~\ref{ass:geometry} and
\ref{ass:initial-region}, for any smooth curve
$\gamma:I\to C_0$, we have
\begin{equation}
\label{eq:hessian-comparison}
\delta_{r_0}\|\dot\gamma(t)\|^2
\le
-\left\langle
\nabla_{\dot\gamma(t)}
\Log_{\gamma(t)}x^\star,
\dot\gamma(t)
\right\rangle_{\gamma(t)}
\le
\zeta_{r_0}\|\dot\gamma(t)\|^2,
\quad t\in I.
\end{equation}
\end{lemma}

\subsection{Definition of the dynamics}
\label{subsec:sigma-ode-definition}

Under Assumptions~\ref{ass:geometry},
\ref{ass:objective}, and \ref{ass:initial-region},
we choose
\begin{equation}
\label{eq:ode-parameters}
\xi\ge\zeta_{r_0},\qquad
0\le\sigma\le\min\{1,\xi-\zeta_{r_0}\},
\qquad
\nu=\frac{\mu}{2\xi-\sigma}.
\end{equation}
Define
\begin{equation}
\label{eq:ode-time-weight}
t_0=
\begin{cases}
\displaystyle
\frac{\operatorname{arsinh}\sqrt{\nu/L}}{\sqrt\nu},
&\nu>0,\\[2mm]
L^{-1/2},&\nu=0,
\end{cases}
\qquad
A(t)=
\begin{cases}
\displaystyle
\frac{\sinh^2(\sqrt\nu(t+t_0))}{\nu},
&\nu>0,\\[2mm]
(t+t_0)^2,&\nu=0.
\end{cases}
\end{equation}
Then
\begin{equation}
\label{eq:ode-weight-identities}
A(0)=\frac1L,\qquad
\dot A=2\sqrt{A(1+\nu A)},\qquad
\ddot A=2+4\nu A.
\end{equation}

We consider the following
\emph{$\sigma$-parameterized Riemannian ITEM ODE}:
\begin{equation}
\label{eq:sigma-item-ode}
\nabla_t\dot X
+
\left(
\frac{(2\xi+1)\dot A}{2A}
+
\frac{(1-\sigma)\nu\dot A}
{2(1+\nu A)}
\right)\dot X
+
2\operatorname{grad}f(X)=0.
\end{equation}
For $\mu>0$, the damping coefficient can be
equivalently expressed as
\begin{equation}
\label{eq:sigma-item-damping}
\sqrt\nu\left[
(2\xi+1)\coth(\sqrt\nu(t+t_0))
+
(1-\sigma)\tanh(\sqrt\nu(t+t_0))
\right].
\end{equation}
When $\mu=0$, the dynamics reduce to
\begin{equation}
\label{eq:sigma-item-convex-ode}
\nabla_t\dot X
+
\frac{2\xi+1}{t+t_0}\dot X
+
2\operatorname{grad}f(X)=0.
\end{equation}
Introducing the rescaled velocity
\begin{equation}
\label{eq:ode-scaled-momentum}
v(t)=\frac{A(t)}{\dot A(t)}\dot X(t),
\end{equation}
the second-order ODE
\eqref{eq:sigma-item-ode}
is equivalent to the first-order system

\begin{equation}
\label{eq:sigma-item-first-order}
\begin{cases}
\displaystyle
\dot X=\frac{\dot A}{A}v,\\[2mm]
\displaystyle
\nabla_t v
=-\frac{\dot A}{A}
\left(
\xi+\frac{(2-\sigma)\nu A}
{2(1+\nu A)}
\right)v
-\frac{\dot A}{2(1+\nu A)}
\operatorname{grad}f(X).
\end{cases}
\end{equation}

We initialize the dynamics by
\begin{equation}
\label{eq:ode-initialization}
X(0)=x_0,\qquad
v(0)=-\frac1L\operatorname{grad}f(x_0).
\end{equation}
The positive time shift ensures that the dynamics
are regular at $t=0$ and that $A(0)=1/L$,
consistent with the subsequent discretization.

\subsection{Connections with existing dynamics}
\label{subsec:sigma-ode-reductions}

The proposed dynamics recover several existing ODEs
under suitable parameter choices and time rescaling.

\paragraph{RNAG-C.}
When $\mu=0$, applying the time transformation
$\tau=\sqrt{2}(t+t_0)$ and
$Y(\tau)=X(\tau/\sqrt{2}-t_0)$ to
\eqref{eq:sigma-item-convex-ode} yields
\begin{equation}
\label{eq:ode-rnag-c-reduction}
\nabla_\tau\dot Y
+\frac{2\xi+1}{\tau}\dot Y
+\operatorname{grad}f(Y)=0.
\end{equation}
This recovers the RNAG-C ODE of
\cite{KimYang2022}.

\paragraph{ITEM and unified AGM.}
In Euclidean space, the formal choice
$(\xi,\sigma)=(1,1)$ gives $\nu=\mu$ and
\begin{equation}
\label{eq:ode-original-item}
\ddot X
+3\sqrt{\mu}
\coth\left(\sqrt{\mu}(t+t_0)\right)\dot X
+2\nabla f(X)=0,
\end{equation}
which recovers the original ITEM ODE
\cite{KimYang2023Kernels,Ushiyama2025}.

Alternatively, choosing $(\xi,\sigma)=(1,0)$
gives $\nu=\mu/2$. Under the same time
transformation, we obtain
\begin{equation}
\label{eq:ode-unified-agm-reduction}
\ddot Y
+\frac{\sqrt{\mu}}{2}
\left[
3\coth\left(\frac{\sqrt{\mu}\tau}{2}\right)
+\tanh\left(\frac{\sqrt{\mu}\tau}{2}\right)
\right]\dot Y
+\nabla f(Y)=0,
\end{equation}
which recovers the unified AGM ODE of
\cite{KimYang2023Unified}.

Note that $(\xi,\sigma)=(1,1)$ lies outside
the admissible range \eqref{eq:ode-parameters},
whereas $(1,0)$ is admissible in the flat case.
These reductions concern the continuous dynamics
and do not imply equivalence of the discrete methods.

\subsection{Lyapunov function and dissipation}
\label{subsec:sigma-ode-lyapunov}

We define the Lyapunov function
\begin{equation}
\label{eq:ode-energy-item}
\begin{aligned}
\mathcal E(t)
={}&A\left(
f(X)-f^\star
-\frac{\sigma\nu}{2}d^2(X,x^\star)
\right)
+(1+\nu A)\left(
\|v-\Log_Xx^\star\|^2
+(\xi-1)d^2(X,x^\star)
\right),
\end{aligned}
\end{equation}
where the time dependence is omitted for brevity. Using \eqref{eq:ode-parameters},
\eqref{eq:ode-weight-identities}, and
\eqref{eq:ode-scaled-momentum}, we obtain
\begin{align}
\mathcal E(t)
&=A(f(X)-f^\star)
+(1+\nu A)\|v-\Log_Xx^\star\|^2
+\left((\xi-1)(1+\nu A)-\frac{\sigma\nu A}{2}\right)d^2(X,x^\star)
\label{eq:ode-energy-squares}\\
&=A(f(X)-f^\star)
+\left(\xi+\frac{\mu A}{2}\right)d^2(X,x^\star)
-\frac{\dot A}{2}\langle\Log_Xx^\star,\dot X\rangle
+\frac{A}{4}\|\dot X\|^2.
\label{eq:ode-energy-expanded}
\end{align}

Since
\[
\xi-1-\frac{\sigma\nu A}{2(1+\nu A)}
\ge\xi-1-\frac{\sigma}{2}
\ge\zeta_{r_0}-1\ge0,
\]
the squared-distance terms in
\eqref{eq:ode-energy-squares} are nonnegative.
Consequently,
\begin{equation}
\label{eq:ode-energy-lower}
\mathcal E(t)
\ge A(t)(f(X(t))-f^\star)
\ge0.
\end{equation}

\begin{lemma}
\label{lem:sigma-ode-dissipation}
Under Assumptions~\ref{ass:geometry},
\ref{ass:objective}, and \ref{ass:initial-region},
let $X(t)$ be a solution of
\eqref{eq:sigma-item-ode} contained in $C_0$.
Then
\begin{equation}
\label{eq:ode-dissipation-identity}
\begin{aligned}
\dot{\mathcal E}
={}&\dot A\left[
f(X)-f^\star
+\langle\operatorname{grad}f(X),\Log_Xx^\star\rangle
+\frac{\mu}{2}d^2(X,x^\star)
\right]\\
&-\frac{\dot A}{2}\left[
\langle\nabla_t\Log_Xx^\star,\dot X\rangle
+\left(\xi+\frac{(1-\sigma)\nu A}{2(1+\nu A)}\right)
\|\dot X\|^2
\right]\le0.
\end{aligned}
\end{equation}
\end{lemma}

\begin{proof}
By the Riemannian differentiation rules,
$\frac{d}{dt}f(X)=\langle\operatorname{grad}f(X),\dot X\rangle$
and
$\frac{d}{dt}d^2(X,x^\star)
=-2\langle\Log_Xx^\star,\dot X\rangle$.
Moreover,
$
\frac{d}{dt}\langle\Log_Xx^\star,\dot X\rangle
=\langle\nabla_t\Log_Xx^\star,\dot X\rangle
+\langle\Log_Xx^\star,\nabla_t\dot X\rangle.
$
Differentiating \eqref{eq:ode-energy-expanded}
term by term gives
\begin{align*}
\dot{\mathcal E}
={}&\dot A\left(f(X)-f^\star+\frac{\mu}{2}d^2(X,x^\star)
+\frac14\|\dot X\|^2\right)
-\left(2\xi+\mu A+\frac{\ddot A}{2}\right)
\langle\Log_Xx^\star,\dot X\rangle\\
&+\left\langle A\operatorname{grad}f(X)
-\frac{\dot A}{2}\nabla_t\Log_Xx^\star,\dot X\right\rangle
+\frac12\left\langle A\dot X-\dot A\Log_Xx^\star,
\nabla_t\dot X\right\rangle.
\end{align*}

Substituting \eqref{eq:sigma-item-ode} and using
$\dot A^2=4A(1+\nu A)$, the gradient--velocity
terms cancel, yielding
\begin{align*}
\dot{\mathcal E}
={}&\dot A\left(f(X)-f^\star
+\langle\operatorname{grad}f(X),\Log_Xx^\star\rangle
+\frac{\mu}{2}d^2(X,x^\star)\right)\\
&-\left(2\xi+\mu A+\frac{\ddot A}{2}
-(2\xi+1)(1+\nu A)-(1-\sigma)\nu A\right)
\langle\Log_Xx^\star,\dot X\rangle\\
&-\frac{\dot A}{2}\left[
\langle\nabla_t\Log_Xx^\star,\dot X\rangle
+\left(\xi+\frac{(1-\sigma)\nu A}{2(1+\nu A)}\right)
\|\dot X\|^2\right].
\end{align*}

Using $\dot A^2=4A(1+\nu A)$ and
$\ddot A=2+4\nu A$ from
\eqref{eq:ode-weight-identities}, the coefficient
of the position--velocity cross term satisfies
\begin{align*}
&2\xi+\mu A+\frac{\ddot A}{2}
-\frac{(2\xi+1)\dot A^2}{4A}
-\frac{(1-\sigma)\nu\dot A^2}{4(1+\nu A)}\\
&\quad=
2\xi+\mu A+1+2\nu A
-(2\xi+1)(1+\nu A)-(1-\sigma)\nu A\\
&\quad=
\left[\mu-(2\xi-\sigma)\nu\right]A=0,
\end{align*}
where the last equality follows from
\eqref{eq:ode-parameters}.
Moreover,
\[
\frac{\dot A}{4}
-\frac{(2\xi+1)\dot A}{4}
-\frac{(1-\sigma)\nu A\dot A}{4(1+\nu A)}
=-\frac{\dot A}{2}
\left(\xi+\frac{(1-\sigma)\nu A}{2(1+\nu A)}\right).
\]
Combining these identities establishes
\eqref{eq:ode-dissipation-identity}.

By \eqref{eq:strong-geodesic-convexity},
\[
f(X)-f^\star
+\langle\operatorname{grad}f(X),\Log_Xx^\star\rangle
+\frac{\mu}{2}d^2(X,x^\star)\le0.
\]
Furthermore,
Lemma~\ref{lem:hessian-comparison-squared-distance}
implies
\[
-\langle\nabla_t\Log_Xx^\star,\dot X\rangle
\le\zeta_{r_0}\|\dot X\|^2.
\]
Consequently,
\[
\dot{\mathcal E}
\le-\frac{\dot A}{2}
\left(
\xi-\zeta_{r_0}
+\frac{(1-\sigma)\nu A}{2(1+\nu A)}
\right)\|\dot X\|^2\le0,
\]
since $\xi\ge\zeta_{r_0}$,
$0\le\sigma\le1$, and $\dot A>0$.
\end{proof}

The additional damping term
$(1-\sigma)\nu\dot A/(2(1+\nu A))$ in
\eqref{eq:sigma-item-ode} ensures the cancellation
of the position--velocity cross term in the
Lyapunov derivative.

\subsection{Region invariance and convergence}
\label{subsec:sigma-ode-rate}

\begin{theorem}
\label{thm:sigma-ode-convergence}
Under Assumptions~\ref{ass:geometry},
\ref{ass:objective}, and \ref{ass:initial-region},
the solution of \eqref{eq:sigma-item-ode} with
initialization \eqref{eq:ode-initialization}
exists uniquely for all $t\ge0$ and satisfies
\begin{equation}
\label{eq:ode-invariance-and-budget}
X(t)\in\overline B_{\mathcal M}(x^\star,r_0),
\qquad
0\le\mathcal E(t)\le\mathcal E(0)\le\xi r_0^2.
\end{equation}
Consequently,
$f(X(t))-f^\star\le\xi r_0^2/A(t)$, yielding

\begin{empheq}[left=\empheqlbrace]{align}
f(X(t))-f^\star
&\le\frac{\xi Lr_0^2}{(1+\sqrt L\,t)^2},
&&\mu=0,
\label{eq:ode-convex-rate}\\
f(X(t))-f^\star
&\le\xi Lr_0^2
\exp\left(-2\sqrt{\frac{\mu}{2\xi-\sigma}}\,t\right),
&&\mu>0.
\label{eq:ode-strongly-convex-rate}
\end{empheq}

\end{theorem}

\begin{proof}
Since $t_0>0$ and $f\in C^2(\mathcal N)$,
\eqref{eq:sigma-item-ode} admits a unique local
solution. If $r_0=0$, the initialized solution
is constant. Hence, assume $r_0>0$.

By geodesic strong convexity,
\[
f(x_0)-f^\star+\frac{\mu}{2}r_0^2
\le
-\langle\operatorname{grad}f(x_0),
\Log_{x_0}x^\star\rangle.
\]

Geodesic convexity and $L$-smoothness imply
$0\preceq\operatorname{Hess}f\preceq LI$
along the minimizing geodesic from $x^\star$
to $x_0$, which yields
\[
\|\operatorname{grad}f(x_0)\|^2
\le-L\langle\operatorname{grad}f(x_0),
\Log_{x_0}x^\star\rangle.
\]

Using $A(0)=1/L$ and
$v(0)=-\operatorname{grad}f(x_0)/L$ in
\eqref{eq:ode-energy-squares}, we obtain
\begin{align*}
\mathcal E(0)-\xi r_0^2
={}&\frac{f(x_0)-f^\star+\frac{\mu}{2}r_0^2}{L}
+\frac{2(1+\nu/L)}{L}
\langle\operatorname{grad}f(x_0),\Log_{x_0}x^\star\rangle+\frac{1+\nu/L}{L^2}
\|\operatorname{grad}f(x_0)\|^2\\
\le{}&\frac{\nu}{L^2}
\langle\operatorname{grad}f(x_0),
\Log_{x_0}x^\star\rangle\le0.
\end{align*}

The preceding initial-energy inequality is strict
unless $\operatorname{grad}f(x_0)=0$.
Indeed, when $\mu>0$, one has $\nu>0$ and
$\langle-\Log_{x_0}x^\star,
\operatorname{grad}f(x_0)\rangle>0$.
When $\mu=0$ and $\operatorname{grad}f(x_0)\ne0$,
convexity along the minimizing geodesic from
$x^\star$ to $x_0$, together with
$\operatorname{grad}f(x^\star)=0$, gives
\[
f(x_0)-f^\star
<-\langle\operatorname{grad}f(x_0),
\Log_{x_0}x^\star\rangle.
\]
Thus, for every nonstationary initialization,

\begin{equation}
\label{eq:ode-strict-initial-energy}
\mathcal E(0)<\xi r_0^2.
\end{equation}

If $\operatorname{grad}f(x_0)=0$, the initial
velocity vanishes and the unique solution is constant;
hence the claimed region invariance is immediate.

For a nonstationary initialization, the trajectory
initially moves into the interior of $C_0$, since
$\dot X(0)=-\dot A(0)\operatorname{grad}f(x_0)$
and
\[
\left.\frac{d}{dt}d^2(X(t),x^\star)\right|_{t=0}
=2\dot A(0)
\langle\operatorname{grad}f(x_0),
\Log_{x_0}x^\star\rangle<0.
\]

Suppose the trajectory subsequently first exits $C_0$.
At its first outward crossing of the boundary,
at some $T>0$,
$d(X(T),x^\star)=r_0$ and
\[
-\langle\Log_{X(T)}x^\star,\dot X(T)\rangle
=\frac12\left.\frac{d}{dt}d^2(X(t),x^\star)
\right|_{t=T}\ge0.
\]
The local comparison lemma and Lyapunov dissipation
hold up to $T$, so by
\eqref{eq:ode-strict-initial-energy},
$\mathcal E(T)\le\mathcal E(0)<\xi r_0^2$.
On the other hand, geodesic strong convexity and
\eqref{eq:ode-energy-expanded} imply
\[
\mathcal E(T)\ge
(\xi+\mu A(T))r_0^2\ge\xi r_0^2,
\]
a contradiction. This argument involves only the
curvature bounds on $C_0$.

To establish global existence, taking the
inner product of \eqref{eq:sigma-item-ode}
with $\dot X$ gives
\[
\frac{d}{dt}
\left(f(X)-f^\star+\frac14\|\dot X\|^2\right)
=
-\frac12
\left(
\frac{(2\xi+1)\dot A}{2A}
+\frac{(1-\sigma)\nu\dot A}{2(1+\nu A)}
\right)\|\dot X\|^2\le0.
\]
Thus $\|\dot X\|$ is bounded. Since $X(t)$
remains in the compact initial ball, the
solution can be extended to all $t\ge0$,
even without completeness of $\mathcal M$.

Finally, \eqref{eq:ode-energy-lower} and
\eqref{eq:ode-invariance-and-budget} yield
\[
f(X(t))-f^\star
\le\frac{\mathcal E(t)}{A(t)}
\le\frac{\xi r_0^2}{A(t)}.
\]
For $\mu=0$, substituting
$A(t)=(t+L^{-1/2})^2$ gives
\eqref{eq:ode-convex-rate}.
For $\mu>0$, the identity
$\dot A/A=2\sqrt{\nu+1/A}\ge2\sqrt\nu$
implies $A(t)\ge L^{-1}e^{2\sqrt\nu t}$,
which proves \eqref{eq:ode-strongly-convex-rate}.
\end{proof}

\section{A Discretization of the Riemannian ITEM Dynamics}
\label{sec:sigma-item-discretization}

We now construct a discrete counterpart of
\eqref{eq:sigma-item-ode} by preserving the structure
of its Lyapunov function. Rather than directly
discretizing the differential equation, we determine
the update coefficients by canceling the mixed terms
in the discrete energy difference.

\subsection{Discrete energy and update rules}
\label{subsec:discrete-construction}

For a fixed step size $s=c/L>0$, initialize
$A_0=1/L$ and
$v_0=-\operatorname{grad}f(x_0)/L$.
Following \eqref{eq:ode-weight-identities} and
\eqref{eq:ode-scaled-momentum}, define
\begin{equation}
\label{eq:discrete-weights}
\begin{aligned}
a_k&=\sqrt{2sA_kB_k},
& A_{k+1}&=A_k+a_k,\\
B_k&=1+\nu A_k,
& \chi_k&=\xi-\frac{\sigma\nu A_k}{2B_k}.
\end{aligned}
\end{equation}
Set
\begin{equation}
\label{eq:discrete-coefficients}
\beta_k=\frac{B_k}{B_{k+1}},\qquad
\eta_k=\frac{a_k}{2B_{k+1}},\qquad
\tau_k=\frac{\beta_k a_k}
{A_k+\chi_k\beta_k a_k}.
\end{equation}

The resulting Riemannian accelerated gradient method
is summarized in Algorithm~\ref{alg:sigma-item}.

\begin{algorithm}[t]
\caption{Region-Preserving Riemannian $\sigma$-ITEM}
\label{alg:sigma-item}
\begin{algorithmic}

\Require Initial point $x_0\in\mathcal N$,
$L>0$, $\mu\in[0,L]$, curvature constants
$\delta_{r_0},\zeta_{r_0}$, admissible parameters $\xi,\sigma$
satisfying \eqref{eq:ode-parameters}, step parameter
$c\in(0,\min\{1,2/\zeta_{r_0}\}]$
satisfying \eqref{eq:discrete-uniform-certificate},
and iteration count $K$.

\State Set
$\displaystyle
\nu=\frac{\mu}{2\xi-\sigma},
\quad s=\frac{c}{L}$.
\State Initialize
$\displaystyle
A_0=\frac1L,\quad
v_0=-\frac1L\operatorname{grad}f(x_0)$.

\For{$k=0,1,\ldots,K-1$}
    \State $\displaystyle
    B_k=1+\nu A_k,\quad
    \chi_k=\xi-\frac{\sigma\nu A_k}{2B_k}$.

    \State $\displaystyle
    a_k=\sqrt{2sA_kB_k},\quad
    A_{k+1}=A_k+a_k$.

    \State $\displaystyle
    B_{k+1}=1+\nu A_{k+1},\quad
    \chi_{k+1}=
    \xi-\frac{\sigma\nu A_{k+1}}{2B_{k+1}}$.

    \State $\displaystyle
    \beta_k=\frac{B_k}{B_{k+1}},\quad
    \eta_k=\frac{a_k}{2B_{k+1}},\quad
    \tau_k=\frac{\beta_k a_k}
    {A_k+\chi_k\beta_k a_k}$.

    \State $y_k=\Exp_{x_k}(\tau_k v_k)$.

    \State $p_k=\mathcal P_{x_k\to y_k}
    \bigl[(1-\chi_k\tau_k)v_k\bigr]$.

    \State $g_k=\operatorname{grad}f(y_k)$.

    \State $x_{k+1}=\Exp_{y_k}(-s g_k)$.

    \State $v_{k+1}=
    \mathcal P_{y_k\to x_{k+1}}
    \bigl[\beta_kp_k-(\eta_k-\chi_{k+1}s)g_k\bigr]$.
\EndFor

\State \textbf{Output:} $x_K$.
\end{algorithmic}
\end{algorithm}

The coefficients in
\eqref{eq:discrete-coefficients} satisfy
\begin{equation}
\label{eq:discrete-cancellation-identities}
\begin{aligned}
2B_{k+1}\eta_k&=a_k,
& B_{k+1}\beta_k&=B_k,\\
A_k\frac{\tau_k}{1-\chi_k\tau_k}
&=\beta_k a_k,
& B_{k+1}\chi_{k+1}-B_k\chi_k
&=\frac{\mu a_k}{2}.
\end{aligned}
\end{equation}

Analogously to \eqref{eq:ode-energy-squares},
we define the discrete Lyapunov function
\begin{equation}
\label{eq:discrete-energy}
\mathcal E_k
=A_k(f(x_k)-f^\star)
+B_k\left(
\|v_k-\Log_{x_k}x^\star\|^2
+(\chi_k-1)d^2(x_k,x^\star)
\right).
\end{equation}
Since $\chi_k\ge\zeta_{r_0}\ge1$, it follows that
\begin{equation}
\label{eq:discrete-energy-lower}
\mathcal E_k\ge A_k(f(x_k)-f^\star)\ge0.
\end{equation}

These identities will cancel the mixed terms
in the discrete energy difference.

\subsection{Discrete Lyapunov analysis}
\label{subsec:discrete-lyapunov}

For the analysis, write
\begin{equation}
\label{eq:discrete-q}
Q_\chi(x,v)
=\chi d^2(x,x^\star)
-2\langle\Log_xx^\star,v\rangle+\|v\|^2,
\end{equation}
so that
$\mathcal E_k=A_k(f(x_k)-f^\star)
+B_kQ_{\chi_k}(x_k,v_k)$.
Let
\[
q_k=\frac{\tau_k}{1-\chi_k\tau_k}
=\frac{\beta_k a_k}{A_k},
\qquad
K_k=\chi_{k+1}-\delta_{r_0}.
\]
In particular,
$\Log_{y_k}x_k=-q_kp_k$.

\begin{lemma}
\label{lem:discrete-energy-descent}
Suppose both geodesic segments in
Algorithm~\ref{alg:sigma-item} lie in $C_0$
and $0<s\le2/(\zeta_{r_0} L)$. Define
\begin{equation}
\label{eq:discrete-descent-coefficients}
\begin{aligned}
P_k={}&B_k(\chi_k-\zeta_{r_0})
(2q_k+\chi_kq_k^2)
+\beta_k\nu a_k
+\frac{\mu A_kq_k^2}{2},\\
\widehat D_k={}&
A_{k+1}s\left(1-\frac{Ls}{2}\right)
-B_{k+1}\eta_k^2-\frac{B_{k+1}K_k}{\chi_{k+1}}
(\chi_{k+1}s-\eta_k)_+^2,
\end{aligned}
\end{equation}
where $(a)_+=\max\{a,0\}$. Then
\begin{equation}
\label{eq:discrete-energy-difference}
\begin{aligned}
\mathcal E_{k+1}-\mathcal E_k
\le{}&
-P_k\|p_k\|^2
+2sB_kK_k\|p_k\|\|g_k\|-\widehat D_k\|g_k\|^2.
\end{aligned}
\end{equation}
Consequently, $\mathcal E_{k+1}\le\mathcal E_k$
whenever
\begin{equation}
\label{eq:discrete-matrix-condition}
P_k\ge0,\qquad
\widehat D_k\ge0,\qquad
P_k\widehat D_k\ge s^2B_k^2K_k^2.
\end{equation}
\end{lemma}

\begin{proof}
By geodesic strong convexity and $L$-smoothness,
\begin{align*}
&A_{k+1}(f(x_{k+1})-f^\star)
-A_k(f(x_k)-f^\star)\\
&\le A_kq_k\langle g_k,p_k\rangle
-a_k\langle g_k,\Log_{y_k}x^\star\rangle
-\frac{\mu a_k}{2}d^2(y_k,x^\star)\\
&\quad-\frac{\mu A_kq_k^2}{2}\|p_k\|^2
-A_{k+1}s\left(1-\frac{Ls}{2}\right)
\|g_k\|^2.
\end{align*}
On the other hand, direct expansion of
\eqref{eq:discrete-q}, together with
\eqref{eq:discrete-cancellation-identities}, gives
\begin{align*}
&B_{k+1}Q_{\chi_{k+1}}
(y_k,\beta_kp_k-\eta_kg_k)
-B_kQ_{\chi_k}(y_k,p_k)\\
&=\frac{\mu a_k}{2}d^2(y_k,x^\star)
+a_k\langle g_k,\Log_{y_k}x^\star\rangle\\
&\quad-A_kq_k\langle g_k,p_k\rangle
-\beta_k\nu a_k\|p_k\|^2
+B_{k+1}\eta_k^2\|g_k\|^2.
\end{align*}
Hence, the distance-square and mixed inner-product
terms cancel exactly.

It remains to bound the change in $Q_\chi$
along the two geodesic updates.
The upper distance-Hessian bound in
Lemma~\ref{lem:hessian-comparison-squared-distance}
yields
\begin{align*}
&Q_{\chi_k}(y_k,p_k)-Q_{\chi_k}(x_k,v_k)\\
&\quad\le-(\chi_k-\zeta_{r_0})
(2q_k+\chi_kq_k^2)\|p_k\|^2.
\end{align*}
Similarly, applying both comparison bounds
along the gradient update gives
\begin{align*}
&Q_{\chi_{k+1}}(x_{k+1},v_{k+1})
-Q_{\chi_{k+1}}
(y_k,\beta_kp_k-\eta_kg_k)\\
&\quad\le
2s\beta_kK_k\|p_k\|\|g_k\|
+\frac{K_k}{\chi_{k+1}}
(\chi_{k+1}s-\eta_k)_+^2\|g_k\|^2.
\end{align*}
Combining these inequalities proves
\eqref{eq:discrete-energy-difference}.
Finally, the quadratic form in
$\|p_k\|$ and $\|g_k\|$ is nonpositive
under \eqref{eq:discrete-matrix-condition}.
\end{proof}

\subsection{Uniform fixed-step certification}
\label{subsec:discrete-step}

The conditions in \eqref{eq:discrete-matrix-condition}
depend on $A_k$. We next obtain a sufficient
condition that is independent of $k$.

Let $s=c/L$, where
$0<c\le c_{\max}:=\min\{1,2/\zeta_{r_0}\}$,
and define
\begin{equation}
\label{eq:discrete-scalar-range}
\alpha=\frac{2c\mu}{(2\xi-\sigma)L},
\qquad
h=\frac{a_k}{A_k}
\in I_c:=[\sqrt{\alpha},\sqrt{\alpha+2c}].
\end{equation}
For $h\in I_c$, set
\begin{equation}
\label{eq:discrete-scalar-functions}
\begin{aligned}
\beta(h)&=\frac{h}{h+\alpha},
&\chi(h)&=\xi-\frac{\sigma\alpha}{2h^2},\\
\chi_+(h)&=\xi-
\frac{\sigma\alpha(1+h)}{2h(h+\alpha)},
&K(h)&=\chi_+(h)-\delta_{r_0}.
\end{aligned}
\end{equation}
When $\mu=0$, we set $\alpha=0$ and use
$\beta=1$ and $\chi=\chi_+=\xi$,
including at $h=0$.

Define
\begin{equation}
\label{eq:discrete-scalar-certificate}
\begin{aligned}
\mathsf P_c(h)
={}&(\chi-\zeta_{r_0})(2\beta h+\chi\beta^2h^2)
+\frac{\alpha}{h+\alpha}+\frac{(2\xi-\sigma)\alpha h^2}
{2(h+\alpha)^2},\\
\mathsf D_c(h)
={}&(1+h)\left(1-\frac c2\right)
-\frac{\beta}{2}-\frac{K}{2\beta\chi_+}
(\chi_+h-\beta)_+^2,
\end{aligned}
\end{equation}
where $\beta,\chi,\chi_+$, and $K$ are
evaluated at $h$. 
Using the identities
\[
\frac{P_k}{B_k}=\mathsf P_c(h),\qquad
\frac{\widehat D_k}{sA_k}=\mathsf D_c(h),\qquad
\frac{2sB_k}{A_k}=h^2,
\]
we obtain the following uniform sufficient condition
for \eqref{eq:discrete-matrix-condition}:
\begin{equation}
\label{eq:discrete-uniform-certificate}
\begin{pmatrix}
\mathsf P_c(h) & -hK(h)/\sqrt{2}\\
-hK(h)/\sqrt{2} & \mathsf D_c(h)
\end{pmatrix}
\succeq 0,
\qquad \forall h\in I_c.
\end{equation}

Thus, a single scalar certificate guarantees
the discrete Lyapunov descent for every iteration.

\subsection{Existence of a certified fixed step size}
\label{subsec:discrete-step-existence}

\begin{proposition}[Existence of a certified step size]
\label{prop:discrete-step-existence}
The uniform step-size condition
\eqref{eq:discrete-uniform-certificate}
admits a positive solution in the following cases:
\begin{enumerate}
\item[(i)] \emph{Flat case.}
If $\delta_{r_0}=\zeta_{r_0}=1$, choose
$\xi=1$ and $\sigma=0$. Then $c=1$ is admissible.

\item[(ii)] \emph{Nonflat case.}
If $\zeta_{r_0}>\delta_{r_0}$, choose
$\xi>\zeta_{r_0}$ and
$0\le\sigma\le\min\{1,\xi-\zeta_{r_0}\}$.
Then any $c\in(0,c_{\mathrm{safe}}]$ is admissible,
where
\begin{equation}
\label{eq:discrete-safe-step}
c_{\mathrm{safe}}
=\min\left\{
\frac14,\frac{2}{\zeta_{r_0}},
\frac{1}{12\xi^2},
\frac{(\xi-\zeta_{r_0}-\sigma/2)^2}
{27(\xi-\delta_{r_0})^4}
\right\}>0.
\end{equation}
\end{enumerate}
\end{proposition}

\begin{proof}
In case (i), $\chi(h)=\chi_+(h)=1$ and
$K(h)=0$. For $c=1$, we have
$\mathsf P_1(h)\ge0$ and
$\mathsf D_1(h)=(1+h-\beta(h))/2\ge0$,
since $0<\beta(h)\le1$. Thus
\eqref{eq:discrete-uniform-certificate} holds.

For case (ii), let
$d=\xi-\zeta_{r_0}-\sigma/2>0$ and
$M=\xi-\delta_{r_0}>0$.
Since $\nu\le L/2$, we have $\alpha\le c$.
For any $0<c\le c_{\mathrm{safe}}$ and
$h\in I_c$, it follows that
\[
h\le\sqrt{3c}\le\frac{d}{3M^2},
\qquad
\beta(h)\ge\frac23,\qquad
\chi_+(h)h\le\xi\sqrt{3c}\le\frac12.
\]
Hence $\chi_+(h)h<\beta(h)$, and the
positive-part term in $\mathsf D_c(h)$ vanishes.
Moreover, $\chi(h)-\zeta_{r_0}\ge d$, yielding
\[
\mathsf P_c(h)\ge\frac{4d}{3}h,\qquad
\mathsf D_c(h)\ge\frac14,\qquad K(h)\le M.
\]
Consequently,
\[
2\mathsf P_c(h)\mathsf D_c(h)
\ge\frac{2d}{3}h
\ge h^2M^2
\ge h^2K(h)^2.
\]
Therefore, \eqref{eq:discrete-uniform-certificate}
holds for all $h\in I_c$.
\end{proof}

Although $c_{\mathrm{safe}}$ provides an explicit
admissible step parameter, it may be conservative.
We therefore define the largest certified step
parameter as
\begin{equation}
\label{eq:discrete-certified-step}
c_{\mathrm{cert}}
:=
\max\left\{
c\in(0,c_{\max}]:
\eqref{eq:discrete-uniform-certificate}
\text{ holds for all }h\in I_c
\right\}.
\end{equation}
For the parameter choices in
Proposition~\ref{prop:discrete-step-existence},
the feasible set is nonempty. Continuity of the
certification matrix and the interval endpoints
ensures that the maximum is attained. The fixed
step $s=c_{\mathrm{cert}}/L$ can be selected
before the iterations.

\subsection{Region invariance and convergence}
\label{subsec:discrete-convergence}

\begin{theorem}[Region invariance and convergence]
\label{thm:discrete-convergence}
Suppose Assumptions~\ref{ass:geometry},
\ref{ass:objective}, and \ref{ass:initial-region}
hold. Let $\xi,\sigma$ satisfy
\eqref{eq:ode-parameters}, and choose
$s=c/L$ with
$c\in(0,\min\{1,2/\zeta_{r_0}\}]$
satisfying \eqref{eq:discrete-uniform-certificate}.

Then Algorithm~\ref{alg:sigma-item} is well-defined
for all $k\ge0$. The iterates $x_k,y_k$ and
both geodesic update segments
$x_k\to y_k$ and $y_k\to x_{k+1}$
remain in $\overline B_{\mathcal M}(x^\star,r_0)$.
Moreover,
\begin{equation}
\label{eq:discrete-region-invariance}
0\le\mathcal E_{k+1}\le\mathcal E_k
\le\mathcal E_0\le\xi r_0^2.
\end{equation}
Consequently,
\begin{equation}
\label{eq:discrete-convergence-rates}
f(x_k)-f^\star\le
\begin{cases}
\displaystyle
\frac{\xi Lr_0^2}
{[1+k(\sqrt{1+\sqrt{2c}}-1)]^2},
&\mu=0,\\[3mm]
\displaystyle
\xi Lr_0^2
\left(
1+\sqrt{\frac{2c\mu}{(2\xi-\sigma)L}}
\right)^{-k},
&\mu>0.
\end{cases}
\end{equation}
\end{theorem}

\begin{proof}
By the initialization in
Algorithm~\ref{alg:sigma-item},
the discrete and continuous Lyapunov functions
have the same initial value,
$\mathcal E_0=\mathcal E(0)$.
Hence Theorem~\ref{thm:sigma-ode-convergence}
gives $\mathcal E_0\le\xi r_0^2$; moreover,
\eqref{eq:ode-strict-initial-energy} yields
$\mathcal E_0<\xi r_0^2$ if
$\operatorname{grad}f(x_0)\ne0$.
If $\operatorname{grad}f(x_0)=0$, the algorithm
is stationary and the claims follow directly.
We henceforth consider the nonstationary case.

We establish region invariance and energy
descent by induction. Suppose
$x_k\in\overline B_{\mathcal M}(x^\star,r_0)$
and $\mathcal E_k\le\xi r_0^2$.
Consider the momentum geodesic
\[
\gamma(\theta)=\Exp_{x_k}(\theta\tau_kv_k),
\qquad 0\le\theta\le1,
\]
with the rescaled momentum
\[
V(\theta)=\mathcal P_{x_k\to\gamma(\theta)}
[(1-\chi_k\theta\tau_k)v_k].
\]
Since $\chi_k\tau_k<1$, $V(\theta)$ and
$\dot\gamma(\theta)$ are positively proportional.
By Lemma~\ref{lem:hessian-comparison-squared-distance},
\[
\begin{aligned}
&\frac{d}{d\theta}
Q_{\chi_k}(\gamma(\theta),V(\theta))\\
&\quad=-2\left\langle
\nabla_\theta\Log_{\gamma(\theta)}x^\star
+\chi_k\dot\gamma(\theta),V(\theta)
\right\rangle\le0,
\end{aligned}
\]
as long as the segment remains in $C_0$,
since $\chi_k\ge\zeta_{r_0}$.

Suppose this segment first exits $C_0$.
At its first outward crossing of the boundary
$d(\gamma(\theta),x^\star)=r_0$, the proportionality
of $V(\theta)$ and $\dot\gamma(\theta)$ implies
$\langle-\Log_{\gamma(\theta)}x^\star,
V(\theta)\rangle\ge0$. Hence
\[
\begin{aligned}
\xi r_0^2
&>\mathcal E_k
\ge B_kQ_{\chi_k}(\gamma(\theta),V(\theta))\\
&\ge B_k\chi_k r_0^2
=\left(\xi+\frac{\mu A_k}{2}\right)r_0^2
\ge\xi r_0^2,
\end{aligned}
\]
a contradiction. Thus the entire momentum
segment remains in $C_0$.

Next, consider the gradient geodesic
$\Exp_{y_k}(-bg_k)$ for $0\le b\le s$.
Geodesic convexity and $L$-smoothness imply
\[
\|g_k\|^2
\le-L\langle g_k,\Log_{y_k}x^\star\rangle.
\]
Applying Lemma~\ref{lem:hessian-comparison-squared-distance}
along each locally defined prefix contained in $C_0$ gives
\begin{align*}
d^2(\Exp_{y_k}(-bg_k),x^\star)
&\le d^2(y_k,x^\star)\\
&\quad-b\left(\frac2L-\zeta_{r_0} b\right)
\|g_k\|^2
\le r_0^2,
\end{align*}
where the last inequality follows from
$s\le2/(\zeta_{r_0} L)$.
If $g_k
e0$, the bound is strict for $0<b<s$,
excluding a first exit before the full step; if
$g_k=0$, the segment is constant. Thus the entire
gradient segment remains in $C_0$.
Assumption~\ref{ass:initial-region}
provides the compactness required to
extend both geodesic updates to their
full prescribed lengths.

Since both update segments lie in $C_0$,
Lemma~\ref{lem:discrete-energy-descent}
and \eqref{eq:discrete-uniform-certificate}
imply $\mathcal E_{k+1}\le\mathcal E_k$,
closing the induction and proving
\eqref{eq:discrete-region-invariance}.

Finally, \eqref{eq:discrete-energy-lower}
and \eqref{eq:discrete-region-invariance} yield
\[
f(x_k)-f^\star
\le\frac{\mathcal E_k}{A_k}
\le\frac{\xi r_0^2}{A_k}.
\]
When $\mu=0$, the recursion
$A_{k+1}=A_k+\sqrt{2sA_k}$ implies
\[
\sqrt{LA_k}
\ge1+k(\sqrt{1+\sqrt{2c}}-1).
\]
For $\mu>0$, we have
\[
\frac{A_{k+1}}{A_k}
=1+\sqrt{2s\left(\frac1{A_k}+\nu\right)}
\ge1+\sqrt{\frac{2c\mu}{(2\xi-\sigma)L}}.
\]
Combining these bounds with $A_0=1/L$
establishes \eqref{eq:discrete-convergence-rates}.
\end{proof}

\section{Conclusions}
\label{sec:conclusions}

In this paper, we developed a region-preserving
accelerated gradient method for smooth geodesically
convex optimization on Riemannian manifolds with
locally bounded sectional curvature.
We first introduced a parameterized Riemannian
inertial dynamical system and established its
Lyapunov dissipation and accelerated convergence.
Guided by the continuous energy structure,
we constructed an explicit discrete method by
eliminating indefinite cross terms in the
discrete Lyapunov difference and derived
verifiable conditions for a fixed step size.

\section*{Acknowledgments}

The author used OpenAI's ChatGPT and Codex
to assist with manuscript drafting and language
editing, literature searches, mathematical
derivations and proof development, and the
implementation and analysis of numerical
experiments. The author assumes full
responsibility for all content.


\bibliographystyle{spmpsci}
\bibliography{references}

\end{document}